\documentclass[12pt,twoside]{article}
\usepackage{amsmath, amssymb, amsthm}
\newtheorem{theorem}{Theorem}

\newtheorem{lemma}{Lemma}

\font\fourteenb=cmb10 at 14pt
\begin{document}
\begin{center}

\vspace*{-1.0cm}\noindent \
 \textbf{MATHEMATICS AND EDUCATION IN MATHEMATICS, 1995}\\[-0.0mm]\
\ Proceedings of the Twenty Fourth Spring Conference of

the Union of Bulgarian Mathematicians,\\[-0.0mm]
\textit{ Svishtov, April   4 - 7, 1995}\\[-0.0mm]
\font\fourteenb=cmb10 at 12pt \ \

   {\bf \LARGE Estimates of the norms of the Toeplitz operators of  $H^{\infty } $
determined by

rational inner functions
   \\ \ \\ \large Peyo Stoilov}
\end{center}

\

\begin{abstract}
Let  $D{\kern 1pt} $  denote the unit disc and  $T{\kern 1pt} $
the unit circle. Let  $E\subseteq D$  and   $b_{n} (E)$  be the
totality of all finite Blaschke prodtucts  $B_{n} $ with  $n$
zeros belonging to  $E$ . Let
 $$\Omega _{n} (E)=\sup \left\{\left\| \frac{1}{2\pi i} \int _{T}\frac{h(\varsigma )}{B_{n} (\varsigma )(\varsigma -z)} {\kern 1pt} {\kern 1pt}d\varsigma  \right\| _{H^{\infty } } :\left\| h\right\| _{H^{\infty } } \le 1,B_{n} \in b_{n} (E)\right\}$$
and $E_{\xi } =\left\{z:\varepsilon \le \left|z\right|<1,
z/\left|z\right|=\xi \right\}{\kern 1pt} ,{\kern 1pt} {\kern 1pt}
{\kern 1pt} {\kern 1pt} {\kern 1pt} {\kern 1pt} {\kern 1pt} {\kern
1pt} {\kern 1pt} \xi \in T,{\kern 1pt} {\kern 1pt} {\kern 1pt}
{\kern 1pt} {\kern 1pt} {\kern 1pt} {\kern 1pt} {\kern 1pt}
\varepsilon
>0.$
An elementary proof is given of the equalities
 $$\Omega _{n} (D)=\Omega _{n} (E_{\xi } )=1+2n,$$
for all   ${\kern 1pt} {\kern 1pt} {\kern 1pt} {\kern 1pt} {\kern
1pt} {\kern 1pt}{\kern 1pt}{\kern 1pt}{\kern 1pt}{\kern 1pt}
 \xi \in T,{\kern 1pt}{\kern 1pt} \varepsilon >0.$
\end{abstract}
\section{Introduction}
 Let  $D{\kern 1pt} $  denote the unit disc and  $T{\kern 1pt} $
the unit circle. Let  $H^{\infty } $ denote the space of all
functios analytic in  $D{\kern 1pt} $  and such that
 $\left\| f\right\| _{H^{\infty } } =\mathop{\sup }\limits_{z\in D} \left|f(z)\right|<\infty .$
 For  $f\in L^{\infty } (T)$ , we denote by  $T_{f} $  the Toeplitz
operator on  $H^{\infty } $ , defined by
 $$T_{f} h=\int _{T}\frac{\overline{f}(\varsigma )h(\varsigma )}{\varsigma -z} {\kern 1pt} {\kern 1pt}\varsigma  dm (\varsigma  ),{\kern 1pt} {\kern 1pt} {\kern 1pt} {\kern 1pt} {\kern 1pt} {\kern 1pt} {\kern 1pt} {\kern 1pt} {\kern 1pt} h\in H^{\infty } .$$
   Here  $m$  denotes normalised Lebesgue measure on  $T$ .

 Let $E\subseteq D$  and  $b_{n} (E)$  be the totality of all finite
Blaschke prodtucts  $B_{n} $ with  $n$  zeros belonging to $E$.
Let
 $$\Omega _{n} (E)=\sup \left\{\left\| T_{B_{n} } \right\| _{H^{\infty } } :B_{n} \in b_{n} (E)\right\}.$$

In the present paper we shall give an elementary proof of the
equality
 $$\Omega _{n} (D)=1+2n,$$
with more complicated method was proven by Pekarski
$\left[1\right]$.

 Previously in $\left[2\right]$  we proved that
 $1+2n\le \Omega _{n} (D)\le 1+\pi n$.
\section{Main results}
Our main result is based on the following lemmas. Let  $f\in
H^{\infty }$ and
 $$\Lambda (f)=\mathop{\sup }\limits_{\eta \in T} \int _{T}\frac{\left|f(\varsigma \eta )-f(\overline{\varsigma }\eta )\right|}{\left|1-\varsigma \right|} dm(\varsigma  )<\infty .$$

\begin{lemma} If  $f\in H^{\infty } $ , then  $\left\| T_{f}
\right\| _{H^{\infty } } \le \left\| f\right\| _{H^{\infty } }
+\Lambda (f).$
\end{lemma}
\begin{proof}

 $$\left\| T_{f} \right\| _{H^{\infty } } =\sup \left\{\mathop{\lim }\limits_{r\to 1-0} \left|\int _{T}\frac{\overline{f}(\varsigma )h(\varsigma )}{1-\overline{\varsigma }r\eta } dm(\varsigma  )\right|:\eta \in T,{\kern 1pt} {\kern 1pt} {\kern 1pt} {\kern 1pt} {\kern 1pt} {\kern 1pt} {\kern 1pt} \left\| h\right\| _{H^{\infty } } \le 1\right\}=$$

 $$=\sup \left\{\mathop{\lim }\limits_{r\to 1-0} \left|\int _{T}\frac{\overline{f}(\varsigma \eta )h(\varsigma \eta )}{1-r\overline{\varsigma }} dm(\varsigma  )\right|:\eta \in T,{\kern 1pt} {\kern 1pt} {\kern 1pt} {\kern 1pt} {\kern 1pt} {\kern 1pt} {\kern 1pt} \left\| h\right\| _{H^{\infty } } \le 1\right\}\le $$

 $$\le\sup \left\{\mathop{\lim }\limits_{r\to 1-0} \left|\int _{T}\frac{\overline{f}(\varsigma \eta )-\overline{f}(\overline{\varsigma }\eta )}{1-r\overline{\varsigma }} h(\varsigma \eta )dm(\varsigma  )\right|:\eta \in T,{\kern 1pt} {\kern 1pt} {\kern 1pt} {\kern 1pt} {\kern 1pt} {\kern 1pt} {\kern 1pt} \left\| h\right\| _{H^{\infty } } \le 1\right\}+\left\| f\right\| _{H^{\infty } } \le $$

 $$\le \sup \left\{\left|\int _{T}\frac{\left|f(\varsigma \eta )-f(\overline{\varsigma }\eta )\right|}{\left|1-\varsigma \right|} dm(\varsigma  )\right|:\eta \in T\right\}+\left\| f\right\| _{H^{\infty } } =\Lambda (f)+\left\| f\right\| _{H^{\infty } } .$$
 We used, that   $g(z)=\overline{f}(\overline{z}\eta )\in H^{\infty } $   and

$$\displaystyle \left|\int _{T}\frac{\overline{f}(\overline{\varsigma }\eta )h(\varsigma \eta )}{1-r\overline{\varsigma }} dm(\varsigma  )\right|\le \left\| f\right\| _{H^{\infty } } \left\| h\right\| _{H^{\infty } } .$$
\end{proof}
\begin{lemma}If  $$ I(z)=\frac{z-a}{1-z\overline{a}} {\kern 1pt} ,{\kern 1pt}
{\kern 1pt} {\kern 1pt} {\kern 1pt} {\kern 1pt} {\kern 1pt} {\kern
1pt} {\kern 1pt} a\in D, {\kern 1pt} {\kern 1pt}{\kern 1pt} {\kern
1pt}  {\kern 1pt} {\kern 1pt} then  {\kern 1pt} {\kern 1pt} {\kern
1pt} {\kern 1pt}{\kern 1pt} \Lambda (I)\le 2.$$
\end{lemma}

\begin{proof}
 $$\Lambda (I)=\mathop{\sup }\limits_{\eta \in T} \int _{T}\left|\frac{\varsigma \eta -a}{1-\varsigma \eta \overline{a}} -\frac{\overline{\varsigma }\eta -a}{1-\overline{\varsigma }\eta \overline{a}} \right|{\kern 1pt} \frac{{\kern 1pt} {\kern 1pt} {\kern 1pt} {\kern 1pt} dm(\varsigma )}{\left|1-\varsigma \right|}  =$$

 $$=\mathop{\sup }\limits_{\eta \in T} \int _{T}\frac{\left(1-\left|a\right|^{2} \right)\left|\varsigma -\overline{\varsigma }\right|{\kern 1pt} dm(\varsigma )}{\left|1-\varsigma \eta \overline{a}\right|\left|1-\overline{\varsigma }\eta \overline{a}\right|\left|1-\varsigma \right|}  \le 2\mathop{\sup }\limits_{\eta \in T} \int _{T}\frac{1-\left|a\right|^{2} }{\left|1-\varsigma \eta \overline{a}\right|\left|1-\overline{\varsigma }\eta \overline{a}\right|}  {\kern 1pt} {\kern 1pt} {\kern 1pt} dm(\varsigma )\le $$

 $$\le 2\mathop{\sup }\limits_{\eta \in T} \left(\int _{T}\frac{1-\left|a\right|^{2} }{\left|1-\varsigma \eta \overline{a}\right|^{2} }  {\kern 1pt} {\kern 1pt} {\kern 1pt} dm(\varsigma )\right)^{1/2} \left(\int _{T}\frac{1-\left|a\right|^{2} }{\left|1-\overline{\varsigma }\eta \overline{a}\right|^{2} }  {\kern 1pt} {\kern 1pt} {\kern 1pt} dm(\varsigma )\right)^{1/2} =2.$$
\end{proof}
\begin{lemma} If  $I_{k} (z)$ ,  $k=1,{\kern 1pt} \,
2,......,n$  is inner functions  $({\kern 1pt} {\kern 1pt}
\left|I_{k} (\varsigma )\right|=1{\kern 1pt} {\kern 1pt} {\kern
1pt} {\kern 1pt} {\kern 1pt} {\kern 1pt} {\kern 1pt} {\kern 1pt}
a.e{\kern 1pt} {\kern 1pt} {\kern 1pt} {\kern 1pt} {\kern 1pt}
on{\kern 1pt} {\kern 1pt} {\kern 1pt} T{\kern 1pt} {\kern 1pt}
{\kern 1pt} )$ , then $\Lambda (I_{1} I_{2} .....I_{n} )\le
\Lambda (I_{1} )+\Lambda (I_{2} )+.....\Lambda (I_{n} ).$
\end{lemma}
\begin{proof}The proof follows at once from the identity

 $$I_{1} (\varsigma \eta )I_{2} (\varsigma \eta ).....I_{n} (\varsigma \eta )-I_{1} (\overline{\varsigma }\eta )I_{2} (\overline{\varsigma }\eta ).....I_{n} (\overline{\varsigma }\eta )=$$

 $$=I_{1} (\varsigma \eta )I_{2} (\varsigma \eta ).....I_{n} (\varsigma \eta )-I_{1} (\overline{\varsigma }\eta )I_{2} (\varsigma \eta ).....I_{n} (\varsigma \eta )+$$

 $$\begin{array}{l} {+I_{1} (\overline{\varsigma }\eta )I_{2} (\varsigma \eta ).....I_{n} (\varsigma \eta )-I_{1} (\overline{\varsigma }\eta )I_{2} (\overline{\varsigma }\eta ).....I_{n} (\varsigma \eta )+} \\ {......................................................................} \end{array}$$

 $$+I_{1} (\overline{\varsigma }\eta )I_{2} (\overline{\varsigma }\eta ).....I_{n-1} (\overline{\varsigma }\eta )I_{n} (\varsigma \eta )-I_{1} (\overline{\varsigma }\eta )I_{2} (\overline{\varsigma }\eta ).....I_{n-1} (\overline{\varsigma }\eta )I_{n} (\overline{\varsigma }\eta )=$$

 $$=\left(I_{1} (\varsigma \eta )-I_{1} (\overline{\varsigma }\eta )\right)I_{2} (\varsigma \eta ).....I_{n} (\varsigma \eta )+$$

 $$\begin{array}{l} {\left(I_{2} (\varsigma \eta )-I_{2} (\overline{\varsigma }\eta )\right)I_{1} (\overline{\varsigma }\eta ).....I_{n} (\varsigma \eta )+} \\ {..................................................} \end{array}$$

 $$+\left(I_{n} (\varsigma \eta )-I_{n} (\overline{\varsigma }\eta )\right)I_{1} (\overline{\varsigma }\eta )I_{2} (\overline{\varsigma }\eta ).....I_{n-1} (\overline{\varsigma }\eta ).$$
\end{proof}
\begin{lemma} If  $B_{n} \in b_{n} (E)$ , then  $\Lambda (B_{n} )\le 2.$
\end{lemma}
\begin{proof} Let $$B_{n} (z)=\prod \limits _{k=1}^{n}\frac{z-a_{k} }{1-z\overline{a_{k} }} {\kern 1pt}  {\kern 1pt} {\kern 1pt} {\kern 1pt} {\kern 1pt} {\kern 1pt} {\kern 1pt} ,{\kern 1pt} {\kern 1pt} {\kern 1pt} {\kern 1pt} {\kern 1pt} {\kern 1pt} {\kern 1pt} {\kern 1pt} a_{k} \in D.$$
From Lemmas 2 and 3 it follows that
$$\Lambda (B_{n} )\le \sum \limits _{k=1}^{n}\Lambda \left(\frac{z-a_{k} }{1-z\overline{a_{k} }} \right){\kern 1pt}  {\kern 1pt} {\kern 1pt} {\kern 1pt} \le 2n{\kern 1pt} .$$
\end{proof}
\begin{theorem}$\Omega _{n} (D)=\Omega _{n} (E_{\xi }
)=1+2n$ {\kern 1pt} {\kern 1pt} for all {\kern 1pt} {\kern 1pt}
${\kern 1pt} \xi \in T,$ $$where {\kern 1pt} {\kern 1pt} {\kern
1pt} {\kern 1pt} {\kern 1pt} {\kern 1pt}E_{\xi }
=\left\{z:\varepsilon \le \left|z\right|<1,{\kern 1pt} {\kern 1pt}
{\kern 1pt} {\kern 1pt} {\kern 1pt} {\kern 1pt}
z/\left|z\right|=\xi \right\}{\kern 1pt} ,{\kern 1pt} {\kern 1pt}
{\kern 1pt} {\kern 1pt} {\kern 1pt} {\kern 1pt} {\kern 1pt} {\kern
1pt}  \varepsilon >0.$$
\end{theorem}
\begin{proof} From Lemmas 1 and 4 follows that

 $$\begin{array}{l} {\Omega _{n} (E_{\xi } )\le \Omega _{n} (D)=\sup \left\{\left\| T_{B_{n} } \right\| _{H^{\infty } } :B_{n} \in b_{n} (D)\right\}\le } \\ {} \\ {\le 1+\sup \left\{\Lambda (B_{n} ):B_{n} \in b_{n} (D)\right\}\le 1+2n.} \end{array}$$
\ \
   Let ${\kern 1pt} \xi \in T.$  We will show that  $$\Omega _{n} (D)=\Omega _{n} (E_{\xi } )=1+2n.$$
\
   Let  $$x_{k} =(1-q^{k} )\xi ,{\kern 1pt} {\kern 1pt} {\kern 1pt} {\kern 1pt} {\kern 1pt} {\kern 1pt} {\kern 1pt} {\kern 1pt} {\kern 1pt} {\kern 1pt} {\kern 1pt} {\kern 1pt} {\kern 1pt} 0<q<1,{\kern 1pt} {\kern 1pt} {\kern 1pt} \varepsilon <1-q , {\kern 1pt} {\kern 1pt} B_{n} (z)=\prod \limits _{k=1}^{n}\frac{z-x_{k} }{1-z\overline{x_{k} }} {\kern 1pt}  {\kern 1pt}. $$  Let ${\kern 1pt} m>n$   and

       $$y_{k} =B'_{n} (x_{k} )(\left|x_{k} \right|^{2} -1)\xi ,{\kern 1pt} {\kern 1pt} {\kern 1pt} {\kern 1pt} {\kern 1pt} {\kern 1pt} {\kern 1pt} {\kern 1pt} k\ne m,{\kern 1pt} {\kern 1pt} {\kern 1pt} {\kern 1pt} {\kern 1pt} {\kern 1pt} {\kern 1pt}y_{m} =B_{n} (x_{m} ).$$
    Since
      $$\left|B'_{n} (z)\right|(1-\left|z\right|^{2} )\le 1,{\kern 1pt} {\kern 1pt} {\kern 1pt} {\kern 1pt} {\kern 1pt} {\kern 1pt} \left|B_{n} (z)\right|\le 1,{\kern 1pt} {\kern 1pt} {\kern 1pt} {\kern 1pt} {\kern 1pt} {\kern 1pt} {\kern 1pt} {\kern 1pt} {\kern 1pt} {\kern 1pt} {\kern 1pt} {\kern 1pt} (z\in D){\kern 1pt} {\kern 1pt} ,$$  then   $\left|y_{k} \right|\le 1$  and
by a well known Carleson interpolation theorem there exists a
function  $h_{0} \in H^{\infty } $ , such that
 $$h_{0} (x_{k} )=y_{k} =B'_{n} (x_{k} )(\left|x_{k} \right|^{2} -1)\xi ,{\kern 1pt} {\kern 1pt} {\kern 1pt} {\kern 1pt} {\kern 1pt} {\kern 1pt} {\kern 1pt} {\kern 1pt} k\ne m,$$
 $$h_{0} (x_{m} )=y_{m} =B_{n} (x_{m} ),$$
 $$\left\| h_{0} \right\| _{H^{\infty } } \le A(q),$$ where   $$A(q)\to 1{\kern 1pt} {\kern 1pt} {\kern 1pt} {\kern 1pt} {\kern 1pt} {\kern 1pt} {\kern 1pt} {\kern 1pt} {\kern 1pt} as{\kern 1pt} {\kern 1pt} {\kern 1pt} {\kern 1pt} {\kern 1pt} {\kern 1pt} q\to 0{\kern 1pt} {\kern 1pt} {\kern 1pt} {\kern 1pt} {\kern 1pt} {\kern 1pt} {\kern 1pt} \left[3\right].$$
Since  $B_{n} \in b_{n} (E_{\xi } ),$ then

 $$\Omega _{n} (E_{\xi } )\ge \left\| T_{B_{n} } \right\| _{H^{\infty } } \ge \frac{1}{A(q)} \frac{1}{2\pi } \left|\int _{T}\frac{h_{0} (\varsigma )}{B_{n} (\varsigma )(\varsigma -x_{m} )} d\varsigma  \right|=$$

 $$=\frac{1}{A(q)} \left|\frac{h_{0} (x_{m} )}{B_{n} (x_{m} )} +\sum \limits _{k=1}^{n}\frac{h_{0} (x_{k} )}{B'_{n} (x_{k} )(x_{k} -x_{m} )}  \right|=\frac{1}{A(q)} \left|1+\sum \limits _{k=1}^{n}\frac{\left|x_{k} \right|^{2} -1}{\left|x_{k} \right|-\left|x_{m} \right|}  \right|.$$
Since  $m>n$ can be every arbitrary long positive integer and
$\left|x_{m} \right|\to 1$  as   $m\to \infty ,$ then

 $$\Omega _{n} (E_{\xi } )\ge \frac{1}{A(q)} \sum \limits _{k=1}^{n}(1+\left|x_{k} \right| )=\frac{1}{A(q)} \left(1+2n-\sum \limits _{k=1}^{n}q^{k}  \right).$$
Using the fact that  $A(q)\to 1{\kern 1pt} {\kern 1pt} {\kern 1pt}
{\kern 1pt} {\kern 1pt} {\kern 1pt} {\kern 1pt} {\kern 1pt} {\kern
1pt} as{\kern 1pt} {\kern 1pt} {\kern 1pt} {\kern 1pt} {\kern 1pt}
{\kern 1pt} q\to 0{\kern 1pt} {\kern 1pt} {\kern 1pt} {\kern 1pt}
{\kern 1pt} {\kern 1pt} ,$   we obtain  $\Omega _{n} (E_{\xi }
)\ge 1+2n.$
This implies
 $\Omega _{n} (D)=\Omega _{n} (E_{\xi } )=1+2n$   for all   ${\kern 1pt} \xi \in T.$

\end{proof}

\

\

\noindent
{\small Department of Mathematics\\
        Technical University\\
        25, Tsanko Dijstabanov,\\
        Plovdiv, Bulgaria\\
        e-mail: peyyyo@mail.bg}

\end{document}